\documentstyle[11pt]{article}
\begin{document}

\title{Stationary Veselov-Novikov
 equation and isothermally asymptotic surfaces
in projective differential geometry}
\author{{\Large Ferapontov E.V.\thanks{
    Present address:
    Fachbereich Mathematik, SFB 288,
    Technische Universit\"at Berlin,
    10623 Berlin,
    Deutschland,\ \ 
    \hbox{e-mail: {\tt fer@sfb288.math.tu-berlin.de}}}}\\
    Institute for Mathematical Modelling\\
    Academy of Science of Russia, Miusskaya 4\\
    125047 Moscow, Russia\\
    e-mail: {\tt fer@landau.ac.ru}}
\date{}
\maketitle

\newtheorem{theorem}{Theorem}

\pagestyle{plain}

\maketitle

\begin{abstract}

It is demonstrated that the stationary Veselov-Novikov
(VN) and the stationary modified Veselov-Novikov (mVN) equations
describe one and the same class of surfaces in projective differential geometry:
the so-called isothermally asymptotic surfaces, examples of which
include arbitrary quadrics and cubics, quartics of Kummer,
projective transforms of affine spheres and rotation  surfaces.
The stationary mVN
equation arises in the Wilczynski approach and plays the role of 
the projective "Gauss-Codazzi" equations,
while the stationary VN equation follows from the 
Lelieuvre representation of surfaces in 3-space. 
This implies an explicit B\"acklund transformation between
the stationary VN and mVN equations which is an analog of the Miura
transformation between their (1+1)-dimensional limits.

\end{abstract}

\newpage

\section{Introduction}

The Veselov-Novikov (VN) equation
\begin{equation}
\begin{array}{c}
u_t=\alpha u_{xxx}+\beta u_{yyy}-3\alpha (vu)_x-3\beta (wu)_y\\
w_x=u_y \\
v_y=u_x \\
\end{array}
\label{VN}
\end{equation}
(here $\alpha, \beta$ are arbitrary constants) was introduced independently in
\cite{Nizhnik}, \cite{VesNov} for real and complex-conjugate $x, y$
respectively and arises from the compatibility conditions of
the linear problem
\begin{equation}
\begin{array}{c}
\nu_{xy}=u\nu \\
\\
\nu_t=\alpha \nu_{xxx}+\beta \nu_{yyy}-3\alpha v\nu_x-3\beta w\nu_y.
\end{array}
\label{lVN}
\end{equation}
The modified Veselov-Novikov (mVN) equation
\begin{equation}
\begin{array}{c}
p_t=\alpha p_{xxx}+\beta p_{yyy}-2\alpha Vp_x-2\beta Wp_y-\alpha pV_x-
\beta pW_y \\
\\
W_x=\frac{3}{2}(p^2)_y \\
\\
V_y=\frac{3}{2}(p^2)_x \\
\end{array}
\label{mVN}
\end{equation}
was introduced in \cite{Bogdanov} and is associated with the two-dimensional
Dirac operator
\begin{equation}
\begin{array}{c}
\psi^1_x=p\psi^2 \\
\psi^2_y=p\psi^1 \\
\\
\psi^1_t=\alpha \psi^1_{xxx}+\beta \psi^1_{yyy}-2\beta W\psi^1_y-
3\alpha p_x\psi^2_x- \beta W_y\psi^1 -2\alpha p V\psi^2 \\
\\
\psi^2_t=\alpha \psi^2_{xxx}+\beta \psi^2_{yyy}-3\beta p_y\psi^1_y-
2\alpha V\psi^2_x- \alpha V_x\psi^2 -2\beta p W\psi^1. \\
\end{array}
\label{lmVN}
\end{equation}
In the (1+1)-dimensional limit $\alpha=\beta=\frac{1}{2}, ~~ y=x$ equations
(\ref{VN}) and (\ref{mVN}) reduce respectively to the KdV and mKdV equations
$$
u_t=u_{xxx}-6uu_x, ~~~~~~ p_t=p_{xxx}-6p^2p_x,
$$
related by the Miura transformation $u=p^2+p_x$.

In this paper we establish a similar link between the stationary VN and 
the stationary mVN
equations, which in the case $\alpha =-\beta$ assume the following forms

\noindent {\bf Stationary VN}:
\begin{equation}
\begin{array}{c}
u_{xxx}-3(vu)_x=u_{yyy}-3(wu)_y\\
w_x=u_y \\
v_y=u_x. \\
\end{array}
\label{sVN}
\end{equation}
\noindent {\bf Stationary mVN}:
\begin{equation}
\begin{array}{c}
p_{xxx}-2Vp_x-pV_x=p_{yyy}-2Wp_y-pW_y \\
\\
W_x=\frac{3}{2}(p^2)_y \\
\\
V_y=\frac{3}{2}(p^2)_x. \\
\end{array}
\label{smVN}
\end{equation}
Let us introduce the linear system
\begin{equation}
\begin{array}{c}
r_{xx}=pr_y+\frac{1}{2}(V-p_y)r \\
\\
r_{yy}=pr_x+\frac{1}{2}(W-p_x)r \\
\end{array}
\label{r}
\end{equation}
the compatibility conditions of which 
coincide with (\ref{smVN}). B\"acklund  transformation between equations
(\ref{sVN}) and (\ref{smVN}) is given by the formulae
\begin{equation}
\begin{array}{c}
u=p^2-2(\ln r_0)_{xy} \\
\\
v=\frac{2}{3}V-2(\ln r_0)_{xx} \\
\\
w=\frac{2}{3}W-2(\ln r_0)_{yy} \\
\end{array}
\label{Backlund}
\end{equation}
where $r_0$ is an arbitrary solution of (\ref{r}). It can be viewed as the
stationary analog of the Miura transformation. We emphasize that although the
(1+1)- and the stationary limits of VN and mVN equations are related by
B\"acklund transformations, there is no explicit link between the full
(2+1)-dimensional equations. 

The origin of this B\"acklund transformation
is purely differential-geometric: both systems (\ref{sVN}) and (\ref{smVN})
arise in projective differential geometry and describe the so-called
isothermally asymptotic surfaces which have been a subject of interest
of differential geometers of the first half of the century. A number of
intriguing properties and interesting examples of isothermally asymptotic
surfaces can be found in the classical textbooks on projective differential
geometry \cite{Fubini}, \cite{Finikov1}, \cite{Finikov2},
\cite{Lane}, \cite{Bol}.

In sect.2 we recall the standard approach to projective differential geometry
following Wilczynski \cite{Wilczynski}
and derive the stationary mVN
equations (\ref{smVN}) playing the role of projective "Gauss-Codazzi" equations
for isothermally asymptotic surfaces. In this approach linear system
(\ref{r}) specifies the radius-vector $r$ of the
surface. Particular examples of isothermally asymptotic surfaces
and exact solutions of the stationary mVN equation corresponding 
to them are discussed in the end of sect.2. These  include quadrics, 
cubics, quartics of Kummer, projective transforms of
surfaces of rotation and affine spheres.

In sect.3 we recall the Lelieuvre representation of surfaces in 3-space and
derive for isothermally asymptotic surfaces the stationary VN equation
(\ref{sVN}). Thus stationary VN and mVN equations are just two different
parametrizations of one and the same class of isothermally asymptotic
surfaces and the corresponding B\"acklund transformation (\ref{Backlund})
follows directly.

\section{Surfaces in projective differential geometry}
                                        
Following \cite{Wilczynski} we define a surface $M^2$
in projective space $P^3$ in terms of solutions of a linear system
\begin{equation}
\begin{array}{c}
r_{xx}=p r_y+\frac{1}{2}(V-p_y) r \\
\\
r_{yy}=q r_x+\frac{1}{2}(W-q_x) r \\
\end{array}
\label{radius}
\end{equation}
where $p, q, V, W$ are certain functions of $x, y$.
Cross-differentiating (\ref{radius}) and assuming
$r, r_x, r_y, r_{xy}$ to be independent, we arrive at the compatibility
conditions
\begin{equation}
\begin{array}{c}
p_{yyy}-2p_yW-pW_y=
q_{xxx}-2q_xV-qV_x\\
\\
W_x=2qp_y+pq_y \\
\\
V_y=2pq_x+qp_x \\
\end{array}
\label{G-C}
\end{equation}
(see also \cite{Lane}, p.120 and \cite{Fer1}). For any fixed solution
$p, q, V, W$ of (\ref{G-C}) the linear system (\ref{radius})  is compatible
and possesses exactly four linearly independent solutions
$r=(r^0, r^1, r^2, r^3)$ which may be regarded as homogeneous coordinates of a
surface in projective space. For definiteness one can think of a surface
$M^2$ in the ordinary 3-space with the radius-vector
$R=(r^1/r^0, r^2/r^0, r^3/r^0)$. Choosing any other four solutions
$\tilde r=(\tilde r^0, \tilde r^1, \tilde r^2, \tilde r^3)$ of the same system
(\ref{radius}) we see that the corresponding surface $\tilde M^2$
with the radius-vector
$\tilde R=(\tilde r^1/\tilde r^0, \tilde r^2/\tilde r^0, \tilde r^3/\tilde r^0)$
is a projective transform of $M^2$, so that any fixed solution
$p, q, V, W$ of equations (\ref{G-C}) defines a surface $M^2$ uniquely up to
projective equivalence. Moreover, a simple calculation yields
\begin{equation}
\begin{array}{c}
R_{xx}=p R_y+ a R_x \\
R_{yy}=q R_x+ b R_y \\
\end{array}
\label{affine}
\end{equation}
$(a=-2 (\ln r^0)_x, ~~ b=-2 (\ln r^0)_y)$ which implies 
that $x, y$ are asymptotic coordinates on $M^2$.
In what follows we assume that our surfaces are
hyperbolic and the corresponding asymptotic coordinates $x, y$ are real
(it is not a problem to reformulate results that follow in the convex
situation regarding $x, y$ as complex conjugate). Since equations (\ref{G-C})
specify a surface uniquely up to projective equivalence, they can be viewed as
the "Gauss-Codazzi" equations in projective geometry.

The most important invariants in projective differential geometry are the
projective metric
\begin{equation}
2pq~dxdy
\label{metric}
\end{equation}
and the Darboux cubic form
\begin{equation}
pdx^3+qdy^3
\label{cubic}
\end{equation}
which define "generic" surface uniquely up to projective equivalence.
Projective metric  (\ref{metric}) gives rise to the projective area functional
\begin{equation}
\int \int ~pq~dxdy.
\label{area}
\end{equation}

Linear system (\ref{affine}) can be the starting point for developing
affine differential geometry of surfaces referred to asymptotic coordinates
$x, y$. Indeed, let $R=(R^1, R^2, R^3)$ be any 3 nonconstant solutions of
(\ref{affine}) viewed as the radius-vector of
a surface $M^2$ in 3-space. Since any other 3 nonconstant solutions
$\tilde R=(\tilde R^1, \tilde R^2, \tilde R^3)$  are related to $R$ through an
affine transformation $\tilde R=AR+b$ where A is a constant $3\times 3$
matrix and $b$ is a constant 3-vector ($R=const$ is always a
solution of (\ref{affine})), system (\ref{affine})  specifies a surface $M^2$
uniquely up to affine equivalence. The compatibility conditions of system
(\ref{affine})  
\begin{equation}
\begin{array}{c}
(q_x+aq+\frac{1}{2}b^2-b_y)_x=2qp_y+pq_y \\
\\
(p_y+bp+\frac{1}{2}a^2-a_x)_y=2pq_x+qp_x \\
\\
a_y=b_x
\end{array}
\label{G-C1}
\end{equation}
 manifest the "Gauss-Codazzi" equations in affine geometry.

{\bf Remark.} The map (\ref{radius})  $\to$ (\ref{affine}) from projective to
affine geometry can be inverted. Let us consider  linear system
(\ref{affine})  satisfying the compatibility conditions (\ref{G-C1}).
In view of $a_y=b_x$ there exists a function $r^0$
satisfying $a=-2(\ln r^0)_x, ~~ b=-2(\ln r^0)_y$. Introducing
the 4-vector $r=(r^0, r^0R^1, r^0R^2, r^0R^3)$ 
we arrive at (\ref{radius})  with $V$ and $W$
given by
\begin{equation}
V=p_y+bp+\frac{1}{2}a^2-a_x, ~~~~~ W=q_x+aq+\frac{1}{2}b^2-b_y.
\label{projinv}
\end{equation}
Thus the "inverse" map (\ref{affine})  $\to$ (\ref{radius})
from affine to projective geometry is just the transformation from affine
invariants $(p, q, a, b)$ to projective invariants $(p, q, V, W)$.

Isothermally asymptotic surfaces in projective differential geometry are
specified by the constraint
\begin{equation}
p=q
\label{p=q}
\end{equation}
in which case linear system (\ref{radius})  assumes the form (\ref{r})
while it's compatibility conditions (\ref{G-C})  reduce to  the
stationary mVN equation (\ref{smVN}). 
The name "isothermally asymptotic" which is due to Fubini \cite{Fubini}
reflects the property that in asymptotic coordinates $x, y$  Darboux's cubic form
(\ref{cubic}) becomes isothermic: $p(dx^3+dy^3)$. These surfaces are projective 
analogs of isothermic surfaces in conformal geometry characterized by the
isothermicity of the metric in coordinates of lines of curvature.
50 years ago the class of isothermally asymptotic surfaces was
probably among the most popular ever discussed in the context of projective
differential geometry. 
Isothermally asymptotic surfaces arise as the focal surfaces of special
W-congruences preserving Darboux's curves 
(zero curves of the Darboux cubic form) and are
uniquely specified by the requirement that  the 3-web formed by
asymptotic and Darboux's curves is hexagonal. We refer to standard
textbooks \cite{Fubini}, \cite{Finikov1}, \cite{Finikov2},
\cite{Lane}, \cite{Bol} for further discussion.

{\bf Remark 1}. Isothermally asymptotic surfaces may be defined by
the equation
\begin{equation}
\left(\ln \frac{p}{q}\right)_{xy}=0
\label{inv}
\end{equation}
which is equivalent to   (\ref{p=q}) in view of 
the form-invariance of system (\ref{radius}) under the following
transformations:
$$
x^{*}=f(x), ~~~ y^{*}=g(y), ~~~ r^{*}=\sqrt {f^{'}g^{'}}~r, ~~~
p^{*}=p{g^{'}}/{(f^{'})^2}, ~~~ q^{*}=q{f^{'}}/{(g^{'})^2},
$$
$$
V^{*}(f^{'})^2=V+S(f), ~~~~ W^{*}(g^{'})^2=W+S(g)
$$
where $S( \cdot )$ is the usual Schwarzian derivative, that is
$$
S(f)=\frac{f{'''}}{f^{'}}-\frac{3}{2}\left(\frac{f^{''}}{f^{'}}\right)^2
$$
(transformation properties of $V$ and $W$ suggest their interpretation as
projective connections along $x-$ and $y-$asymptotic curves, respectively).
So  (\ref{inv}) is just the invariant form of  (\ref{p=q}). 
It is important that fixing the normalization $p=q$ we  fix coordinates 
$x, y$ up to affine transformations
$x^{*}=ax+b, ~ y^{*}=ay+c, ~ a, b, c = const$.
Thus any isothermally asymptotic surface is endowed with a canonical affine
structure. In these coordinates Darboux's curves are just straight lines
$x+y=const$. The lines $x+ky=const$ also make geometric sense: they can be 
defined in an invariant way as those curves on $M^2$ which have constant 
cross-ratio with the asymptotic and Darboux's directions.

{\bf Remark 2}. Equations (\ref{smVN}) are invariant under the discrete symmetry
$(p, V, W)\to (-p, V, W)$. Geometrically, this means that the class
of isothermally asymptotic surfaces is self-dual. Indeed, two surfaces 
(\ref{r}) corresponding to $(p, q, V, W)$ and $(-p, -q, V, W)$
are projective duals of each other.

For isothermally asymptotic surfaces projective metric assumes the form
$2p^2~dxdy$ with the corresponding area functional
$\int \int p^2~dxdy$ which is the conserved quantity of the mVN equation
(\ref{mVN}).

Let us list some of the most important examples of 
isothermally asymptotic surfaces with the emphasize on solutions of
the stationary mVN equation (\ref{smVN}) corresponding to them.

\bigskip

{\bf Quadrics} correspond to the trivial solution $p=0, ~ W=W(y), ~ V=V(x)$.
\bigskip

{\bf Projective transforms of rotation surfaces} are specified by
$p=p(x+y), ~ W=V=\frac{3}{2}p^2+c$ where p is an arbitrary function of $(x+y)$ and
$c$ is an arbitrary constant. For $c>0$ these are indeed projective transforms
 of surfaces $z=f(x^2+y^2)$, while the cases $c=0$ and $c<0$
correspond to projective transforms of 
 surfaces $z=f(x^2+y)$ and $z=f(x^2-y^2)$, respectively.
Travelling-wave solutions $p(x+cy)$
of equation (\ref{smVN}) correspond to surfaces, which are invariant 
under one-parameter groups of projective transformations. In the case $c\ne 1$
the function $p$ is no longer arbitrary and can be expressed in elliptic functions
(compare with \cite{Lingenberg1}).

\bigskip

{\bf Cubic surfaces} are specified by the following additional constraints in 
(\ref{smVN}):
\begin{equation}
V=-\frac{1}{2}(\ln p)_{xx}+\frac{1}{8}(\ln p)_x^2+\frac{5}{2}p_y, ~~~
W=-\frac{1}{2}(\ln p)_{yy}+\frac{1}{8}(\ln p)_y^2+\frac{5}{2}p_x
\label{1}
\end{equation}
\cite{Lane1}, see also \cite{Lane}, p.131. With these $V,W$ equations 
(\ref{smVN}) imply
$$
\left(\frac{(\ln p)_{xy}}{\sqrt p}+4p\sqrt p\right)_y=5\frac{p_{xx}}{\sqrt p},
~~~
\left(\frac{(\ln p)_{xy}}{\sqrt p}+4p\sqrt p\right)_x=5\frac{p_{yy}}{\sqrt p}.
$$
Integration of these equations for $p$ would provide a 4-parameter family
of exact solutions of equation (\ref{smVN}): indeed, up to projective 
equivalence cubics in $P^3$ depend on 4 essential parameters.
\bigskip

{\bf The Roman surface of Steiner} is a rational quartic in $P^3$ with the 
equation
$$
(x^2+y^2+z^2-1)^2=((z-1)^2-2x^2)((z+1)^2-2y^2)
$$
owing it's name to   Steiner who investigated this surface
in Rome in 1844. Besides quadrics and ruled cubic surfaces
the Roman surface of Steiner is the only surface in $P^3$ 
possessing infinitely many conic sections 
through any of it's points. This result was announced
several times: by Moutard in 1865,  Darboux in 1880 and Wilczynski in 1908 (see
\cite{Wilczynski}, 1909 for historical remarks). The interest to the  
Roman surface of Steiner in projective differential geometry is due to 
the remarkable
construction of Darboux, relating with an arbitrary surface $M^2$ in $P^3$
and an arbitrary point $p$ on $M^2$ an osculating Roman surface of Steiner
which has the fourth order of tangency with $M^2$ at this point.
Analytically, the Roman surface of Steiner corresponds to the choice
\begin{equation}
V=-\frac{1}{2}(\ln p)_{xx}+\frac{1}{8}(\ln p)_x^2-\frac{5}{2}p_y, ~~~
W=-\frac{1}{2}(\ln p)_{yy}+\frac{1}{8}(\ln p)_y^2-\frac{5}{2}p_x,
\label{2}
\end{equation}
$$
(\ln p)_{xy}=\frac{4}{9}p^2
$$
(\cite{Bol}, p.149-150) 
implying upon substitution in (\ref{smVN})  the following equations for $p$:
\begin{equation}
\begin{array}{c}
p_{xx}=-\frac{4}{3}pp_y\\
\\
p_{yy}=-\frac{4}{3}pp_x\\
\\
(\ln p)_{xy}=\frac{4}{9}p^2.
\end{array}
\label{Roman}
\end{equation}
These can be explicitely integrated:
$$
p^2=\frac{9}{4}\frac{f^{'}g^{'}}{(f+g)^2}
$$
where the functions $f(x)$ and $g(y)$ satisfy the ODE's
$$
(f^{'})^3=(a_0+a_1f+a_2f^2)^2, ~~~ (g^{'})^3=(a_0-a_1g+a_2g^2)^2.
$$
Here $a_i$ are arbitrary constants. 
Under the transformation $(p, V, W)\to (-p, V, W)$ equations
(\ref{2}) transform to (\ref{1}). This means, that the dual of the Roman 
surface of Steiner is a cubic, and hence the Roman surface itself is a quartic
of class 3 (\cite{Bol}, p.150).

The Roman surface of Steiner belongs to
a broader class of isothermally asymptotic quartic surfaces known as

\bigskip

{\bf Quartics of Kummer} investigated by Kummer as singular surfaces of 
quadratic line complexes. Around 1870 Kummer himself constructed plaster
models of his quartics which now belong to the G\"ottingen collection
(pictures of these models and the necessary explanations
can be found in \cite{Fischer}). Analytically, the quartics of Kummer 
are specified by the conditions
\begin{equation}
\begin{array}{c}
V=\frac{11}{8}(\ln p)_{xx}+2(\ln p)_x^2, ~~~
W=\frac{11}{8}(\ln p)_{yy}+2(\ln p)_y^2, \\
\ \\
(\ln p)_{xy}=\frac{4}{9}p^2
\end{array}
\label{3}
\end{equation}
(\cite{Bol}, p.231). Substituting these $V, W$ in  (\ref{smVN}) we arrive at
\begin{equation}
\left(\frac{1}{p^2}(p^2(p^2)_y)_y\right)_y=
\left(\frac{1}{p^2}(p^2(p^2)_x)_x\right)_x.
\label{p}
\end{equation}
With
$$
p^2=\frac{9}{4}\frac{f^{'}g^{'}}{(f+g)^2}
$$
equations (\ref{p}) can be rewritten in the form
\begin{equation}
\begin{array}{c}
\frac{1}{3}(f+g)^3
\left(\frac{d^3((f^{'})^3)}{df^3}-\frac{d^3((g^{'})^3)}{dg^3}\right)-
4(f+g)^2
\left(\frac{d^2((f^{'})^3)}{df^2}-\frac{d^2((g^{'})^3)}{dg^2}\right)+ \\
\ \\
20(f+g)
\left(\frac{d((f^{'})^3)}{df}-\frac{d((g^{'})^3)}{dg}\right)-
40\left((f^{'})^3-(g^{'})^3\right)=0
\end{array}
\label{dfg}
\end{equation}
(here we used the identities $\partial_x=f^{'}\partial_f, ~ 
\partial_y=g^{'}\partial_g$). Applying to (\ref{dfg})
operator
$\partial^6/ \partial f^3\partial g^3$ we arrive at
$$
\frac{d^6((f^{'})^3)}{df^6}-\frac{d^6((g^{'})^3)}{dg^6}=0
$$
implying that $(f^{'})^3$ and $(g{'})^3$ are polynomials of the 6-th order in 
$f$ and $g$, respectively. Coefficients of these polynomials are not 
independent and can be fixed upon substitution in (\ref{dfg}):
\begin{equation}
(f^{'})^3=P(f), ~~~~
(g^{'})^3=P(-g)
\label{fg}
\end{equation}
where $P$ is an arbitrary polynomial of the 6th order. 
Calculations we have presented here follow \cite{Finikov1}, p.66-69.
Formulae (\ref{fg}) reflect the uniformizability of Kummer's quartics via
theta functions of genus 2 \cite{Hudson}. Since equations (\ref{3}) are invariant
under the transformation $(p, V, W)\to (-p, V, W)$
the class of Kummer's quartics is self-dual.

Quartics of Kummer constitute  a subclass of

\bigskip

{\bf Projectively applicable surfaces}
which are characterized by a condition
\begin{equation}
(\ln p^2)_{xy}=cp^2
\label{c}
\end{equation}
for some constant $c$. Quartics of Kummer correspond to $c=\frac{8}{9}$.
Condition (\ref{c}) means that the projective metric $2p^2~dxdy$ has 
constant Gaussian curvature $K=-c$. To investigate equations (\ref{smVN})
with the additional constraint (\ref{c}) we introduce the anzatz
\begin{equation}
V=(\ln p)_{xx}+\frac{1}{2}(\ln p)_x^2+A, ~~~
W=(\ln p)_{yy}+\frac{1}{2}(\ln p)_y^2+B,
\label{AB}
\end{equation}
which implies upon substitution in (\ref{smVN}) the following equations 
for $A, B$:
\begin{equation}
\begin{array}{c}
A_y=\frac{3}{2}(1-c)(p^2)_x, ~~~~ B_x=\frac{3}{2}(1-c)(p^2)_y, \\
\ \\
(p^2)_xA+p^2A_x=(p^2)_yB+p^2B_y.
\end{array}
\label{dAB}
\end{equation}
For $c=1$ equations (\ref{dAB}) are satisfied if $A=B=0$. The corresponding 
surfaces are improper affine spheres; they will be discussed below. Here we 
consider the case $c\ne 1$. Introducing $F$ by the formulae
\begin{equation}
A_x=-A(\ln p^2)_x+F, ~~~~ B_y=-B(\ln p^2)_y+F
\label{F}
\end{equation}
and writing down the compatibility conditions of (\ref{F}) with
$(\ref{dAB})_1$, $(\ref{dAB})_2$ we arrive at the equations for $F$
\begin{equation}
\begin{array}{c}
F_x=2cBp^2+\frac{3}{2}(1-c)\frac{1}{p^2}(p^2(p^2)_y)_y \\
\ \\
F_y=2cAp^2+\frac{3}{2}(1-c)\frac{1}{p^2}(p^2(p^2)_x)_x \\
\end{array}
\label{dF}
\end{equation}
the compatibility conditions of which coincide with (\ref{p}). Inserting in 
(\ref{p})
\begin{equation}
p^2=\frac{1}{c}\frac{f^{'}g^{'}}{(f+g)^2}
\label{pc}
\end{equation}
we end up with the same $f, g$ as in (\ref{fg}). For any $p$ given
by (\ref{pc}) the functions 
$V$ and $W$ can be recowered from (\ref{AB}) where $A, B, F$
satisfy the compatible system (\ref{dAB}), (\ref{F}), (\ref{dF}). 
Thus for any such $p$ there exists 3-parameter family of surfaces which 
have the same metric $2p^2~dxdy$ and the same cubic form
$p(dx^3+dy^3)$ and are not projectively equivalent. In general, two 
projectively different surfaces in $P^3$ having the same projective metric
(\ref{metric}) and the same cubic form (\ref{cubic}) in a common asymptotic
parametrization $x, y$ are called projectively applicable. One can show 
that only for isothermally asymptotic surfaces $M^2$ satisfying 
(\ref{c}) does there exist 3-parameter family of projectively different
surfaces which are all projectively applicable to $M^2$ 
(the value 3 is the maximal possible). We point out that geometry of 
the surface $M^2$ depends crucially on the value of constant $c$.
Isothermally asymptotic surfaces possessing only one-parameter families of
projective applicabilities have been discussed in \cite{Lingenberg2}.

\bigskip

{\bf Affine spheres} constitute an important subclass of isothermally
asymptotic surfaces specified by the following reduction in (\ref{smVN}):
\begin{equation}
V=\frac{p_{xx}}{p}-\frac{1}{2}\left(\frac{p_x}{p}\right)^2, ~~~~~
W=\frac{p_{yy}}{p}-\frac{1}{2}\left(\frac{p_y}{p}\right)^2.
\label{VW}
\end{equation}
After this ansatz the first equation in (\ref{smVN}) will be satisfied
identically while the last two imply the Tzitzeica equation for $p$:
\begin{equation}
(\ln p)_{xy}=p^2+\frac{c}{p}, ~~~ c=const.
\label{Tzitzeica}
\end{equation}
The cases $c\ne 0$ and $c=0$ correspond to proper and improper affine spheres,
respectively. With $V, W$ given by (\ref{VW}) equations (\ref{r}) possess
particular solution $r^0=\sqrt p$. With $R=r/r^0$ equations (\ref{affine})
assume the form
$$
\begin{array}{c}
R_{xx}=p R_y-\frac{p_x}{p} R_x \\
\ \\
R_{yy}=p R_x-\frac{p_y}{p} R_y \\
\end{array}
$$
which become the familiar equations for the radius-vector of affine spheres
after adding the compatible equation
$$
R_{xy}=-\frac{c}{p}R.
$$
The fact that Tzitzeica's equation solves the stationary mVN equation is also
reflected in the following nonlocal representation of the mVN
equation (in case $\alpha =-\beta =1$):
$$
p_t=\left(\frac{1}{p}\partial_x p^2 \partial_y^{-1}\frac{1}{p}\partial_x -
  \frac{1}{p}\partial_y p^2 \partial_x^{-1}\frac{1}{p}\partial_y\right)
  \left(p(\ln p)_{xy}-p^3\right).
$$
\bigskip

In the recent paper \cite{KonPin}
Konopelchenko and Pinkall introduced integrable dynamics of surfaces in
3-space governed by the VN equation. One can show that
isothermally asymptotic surfaces can be interpreted as the stationary points
of this evolution. So it is not surprising that they are described
by the stationary mVN equation (which, as will be 
demonstrated below, is equivalent to the stationary VN).

Similar class of surfaces, described by the stationary mVN equation (with
different real reduction) arises in Lie sphere geometry \cite{Fer}.
These are the
so-called diagonally cyclidic surfaces which in Lie sphere geometry play a
role similar to that of isothermally asymptotic surfaces in projective
geometry.

\section{The Lelieuvre representation of surfaces in 3-space}

We consider a surface $M^2\in E^3$ parametrized by asymptotic coordinates $x,
y$ with the radius-vector $R(x, y)$ satisfying equations
(\ref{affine})
$$
\begin{array}{c}
R_{xx}=p R_y+ a R_x \\
R_{yy}=q R_x+ b R_y \\
\end{array}
$$
with the compatibility conditions (\ref{G-C1}).

The Lelieuvre representation of a surface $M^2$ \cite{Lelieuvre}
is defined by the formulae
\begin{equation}
\begin{array}{c}
R_x=\nu\times \nu_x \\
R_y=-\nu\times \nu_y \\
\end{array}
\label{Lelieuvre}
\end{equation}
where a 3-vector $\nu$ (called the affine conormal) satisfies the equation
\begin{equation}
\nu _{xy}=u\nu
\label{Schroedinger}
\end{equation}
for certain potential $u(x, y)$. In formulae  (\ref{Lelieuvre})
"$\times$" denotes vector product in $E^3$. Equations (\ref{Lelieuvre})  are
compatible in view of (\ref{Schroedinger}).  
Our aim is to show that in the case  of isothermally asymptotic surfaces 
potential $u$ satisfies the stationary VN equation (\ref{sVN}).

Assuming $\nu, \nu_x, \nu_y$ to
be independent we introduce the expansions
$$
\begin{array}{c}
\nu_{xx}=c^1 \nu_x+ c^2 \nu_y +c^3 \nu \\
\nu_{yy}=c^4 \nu_x+ c^5 \nu_y +c^6 \nu \\
\end{array}
$$
where the coefficients $c^i$ are uniquely expressible through $p, q, a, b$.
Indeed, differentiating the first equation (\ref{Lelieuvre})  with respect
to $x$ and the second with respect to $y$ we arrive at $c^1=a,~ c^2=-p,~ c^4=-q,~
c^5=b$. The remaining coefficients can be obtained from the compatibility
conditions of the equations for $\nu$:
$u=a_y+pq=b_x+pq, ~~ c^3=p_y+pb, ~~ c^6=q_x+qa$,  so that the final equations
for $\nu$ assume the form
\begin{equation}
\begin{array}{c}
\nu _{xy}=u\nu \\
\\
\nu_{xx}=a \nu_x-p\nu_y + (p_y+pb)\nu \\
\\
\nu_{yy}=-q\nu_x+b \nu_y + (q_x+qa)\nu
\end{array}
\label{conormal}
\end{equation}
($u=a_y+pq=b_x+pq$). The compatibility conditions of equations (\ref{conormal})
coincide with (\ref{G-C1}).

{\bf Remark.} Introducing $r^0$ by the formulae $a=-2(\ln r^0)_x, ~
b=-2(\ln r^0)_y $ one can easily check that $\nu=(r^0)^2(R_x\times R_y)$
satisfies  (\ref{conormal}).

In the case $p=q$ of isothermally asymptotic surfaces the radius-vector $R$
satisfies the equations
$$
\begin{array}{c}
R_{xx}=p R_y+ a R_x \\
R_{yy}=p R_x+ b R_y \\
\end{array}
$$
with  the compatibility conditions
$$
\begin{array}{c}
(p_x+ap+\frac{1}{2}b^2-b_y)_x=\frac{3}{2}(p^2)_y \\
\\
(p_y+bp+\frac{1}{2}a^2-a_x)_y=\frac{3}{2}(p^2)_x \\
\\
a_y=b_x
\end{array}
$$
while equations for the conormal $\nu$ assume the form
\begin{equation}
\begin{array}{c}
\nu _{xy}=u\nu \\
\\
\nu_{xx}=a \nu_x-p\nu_y + (p_y+pb)\nu \\
\\
\nu_{yy}=-p\nu_x+b \nu_y + (p_x+pa)\nu
\end{array}
\label{conormal1}
\end{equation}
($u=a_y+p^2=b_x+p^2$). It is a direct check that in this case equations
(\ref{conormal1}) imply
\begin{equation}
\begin{array}{c}
\nu _{xy}=u\nu \\
\nu_{xxx}-3v \nu_x=
\nu_{yyy}-3w\nu_y\\
v_y=u_x\\
w_x=u_y
\end{array}
\label{lsVN}
\end{equation}
where $v=\frac{2}{3}V+a_x, ~~ w=\frac{2}{3}W+b_y$ and $V, W$ are 
projective invariants (\ref{projinv}).
Since (\ref{lsVN}) is just the stationary VN linear problem, 
the stationary VN equation (\ref{sVN}) for
$u, v, w$ follows directly. Thus any solution $p, V, W$ of the stationary mVN equation
generates solution $u, v, w$ of the stationary VN equation according to the
formulae
$$
\begin{array}{c}
u=p^2+a_y=p^2+b_x \\
\\
v=\frac{2}{3}V+a_x\\
\\
w=\frac{2}{3}W+b_y\\
\end{array}
$$
where $a=-2(\ln r^0)_x, ~ b=-2(\ln r^0)_y $ and $r^0$ is a solution of
(\ref{r}). This gives the desired B\"acklund transformation (\ref{Backlund}).

Let us apply transformation (\ref{Backlund}) to affine spheres
choosing $r^0=\sqrt p$ which is a solution of linear system
(\ref{r}). With this $r^0$ formulae (\ref{Backlund}), (\ref{VW}),
(\ref{Tzitzeica}) give the following expressions for $u, v, w$:
$$
u=-\frac{c}{p}, ~~~~~ v=-\frac{1}{3}\frac{p_{xx}}{p} + \frac{2}{3}
\left(\frac{p_x}{p}\right)^2, ~~~~~
w=-\frac{1}{3}\frac{p_{yy}}{p} + \frac{2}{3}
\left(\frac{p_y}{p}\right)^2
$$
implying that
$$
(\ln u)_{xy}=u-\frac{c^2}{u^2}, ~~~~  v=\frac{1}{3}\frac{u_{xx}}{u}, ~~~~
w=\frac{1}{3}\frac{u_{yy}}{u}
$$
solve (\ref{sVN}). 
This is reflected in the  nonlocal representation of the VN
equation (in the case $\alpha =-\beta =1$)
$$
u_t=\left(\partial_x u \partial_y^{-1}\frac{1}{u^2}\partial_x -
  \partial_y u \partial_x^{-1}\frac{1}{u^2}\partial_y\right)
  \left(u^2(\ln u)_{xy}-u^3\right)
$$
which was obtained in \cite{KonRog}.

As it was pointed out in \cite{Matveev}, 
VN equation is covariant under the Moutard
transformation which can be specialized to the stationary case as well.
Geometrically, this construction should give the known B\"acklund
transformation for isothermally asymptotic surfaces which is generated by
a W-congruence preserving Darboux's curves (see e.g. \cite{Finikov2}).

\section{Acknowledgements}

This research was supported by the RFFI grants 96-01-00166,
96-06-80104, INTAS 96-0770 and the Alexander von Humboldt Foundation.
The author is  thankful to B.G. Konopelchenko for useful discussions.

\end{document}